\documentclass{amsart}

\newtheorem{theorem}{Theorem}[section]

\theoremstyle{remark}
\newtheorem{remark}[theorem]{Remark}

\numberwithin{equation}{section}

\newtheorem{proposition}[theorem]{Proposition}

\begin{document}

\date{}
 
\title{Nonlinear wave equations 
and   singular solutions}

\author{Hideshi YAMANE}
\thanks{
  This research was partially supported by Grant-in-Aid for 
  Scientific Research (No.17540182), Japan 
  Society for the Promotion of Science. 
}

\address{Department of Physics,  Kwansei Gakuin 
University, 
Gakuen 2-1, Sanda, Hyougo 669-1337, Japan}

\email{yamane@ksc.kwansei.ac.jp}


\begin{abstract}
We construct solutions to  
nonlinear wave equations that are  
singular along    
a prescribed noncharacteristic hypersurface 
which is  the graph of a function 
satisfying not the Eikonal but another 
partial differential equation of the first order. 
The method of Fuchsian reduction is 
employed. 
\end{abstract}

\maketitle

\section{Introduction}
\label{sec:intro}

In our previous article \cite{Yam}, we constructed 
a solution 
to a nonlinear Cauchy problem and gave 
a lower bound of its lifespan.

In the present article, we   consider certain  
nonlinear wave equations with  polynomial nonlinearity. 
More precisely, the right hand side of an equation in our class 
is a polynomial of  degree $2$ in the derivatives of the unknown 
function. 
Notice that some of this kind of  equations  fall within the 
class treated in \cite{Yam}. 

We  construct solutions 
that are singular along a prescribed 
\textit{noncharacteristic} hypersurface 
satisfying a certain condition written in terms of 
a partial differential equation of the first order. 
Blowup occurs due to the presence 
of a logarithmic term. Its 
coefficient (and hence the 
speed of blowup) is determined by the hypersurface and 
the nonlinear term of the equation. 
On the  other hand, it is possible to  construct 
hypersurfaces in accordance with a given pair of the 
coefficient and the nonlinear term.

In later sections, we  briefly discuss  a nonlinear 
elliptic equation and multi-valued solutions to 
partial differential equations in the complex domain. 
We also consider equations with nonlinearity of 
higher order.

The proofs are based on   Fuchsian reduction. 
This method, a powerful tool in many branches of analysis,    
was first employed 
by Kichenassamy-Littman in their construction 
of singular solutions to the equation  
 $\square u = \exp u$ in \cite{KL1}.  
Similar methods were used in \cite{Kobayashi}, 
\cite{Tahara1} and \cite{Tahara2}. 
We reduce our    equations to 
Fuchsian ones   and use   results by G\'erard-Tahara in \cite{GT_art} 
and \cite{GT_book}(Chapter 8). 
An alternative proof is the reduction to  first-order 
Fuchsian systems which are solved in  \cite{KL1} and \cite{KL2}.

\section{Main result}
\label{sec:result}

Let $f(t, x; \tau, \xi)$ be a polynomial of 
degree $2$ in 
$(\tau, \xi)\in\mathbb{R}\times\mathbb{R}^n$  whose 
coefficients  are analytic functions in 
$(t, x)$ in an open set $U$ of $\mathbb{R}\times\mathbb{R}^n$. 
We decompose $f$ into homogeneous parts: 
$f=f_2+f_1+f_0$, where $f_l$ is   homogeneous of 
degree $l$. For $l=1, 2$, we sometimes denote 
$f_l (t, x; \tau, \xi)$ by $f_l (\tau, \xi)$ for brevity. 
Notice that $f_0=f_0 (t, x)$ is just a function in $(t, x)$ 
and is independent of $(\tau, \xi)$. 

The primary purpose of the present paper is the  study 
of a nonlinear wave equation 
\begin{equation}
 u_{tt}-\Delta u=f(t, x; \partial_t u, \nabla u)
 =\sum_{l=0}^2 f_l (t, x; \partial_t u, \nabla u), 
\label{eq:wave}
\end{equation}
where $\nabla$ is the gradient with respect to $x$.

We   construct a solution 
$u(t, x)$   which  blows up on a prescribed 
noncharacteristic hypersurface. 
More precisely, our main result is the following:
\begin{theorem}
\label{thm:main}
 Let $\Sigma$ be a hypersurface defined by 
 \(\Sigma=\{ 
             (t, x)\in U;     t=\psi(x)
          \}
 \), 
 where  $\psi(x)$ is an analytic function. 
 Assume 
 \begin{equation}
 \label{eq:assumption}
    0\ne 1-|\nabla\psi(x)|^2 = 
    a f_2 \bigl(\psi(x), x; -1, \nabla \psi(x)\bigr) 
 \end{equation}
 in $U$ for  a nonzero constant $a$.  
 Then for any analytic function 
 $v_0$ on $\Sigma$, 
 there exists an  analytic function $v(t, x)$ 
 in an open neighborhood $\Omega$ of $\Sigma$ 
  such that: 
 \begin{itemize}
 \item
   \(
     u(t, x)= -a \log  (t-\psi(x)) + v(t, x)
   \) 
   is a solution to  
   (\ref{eq:wave}) in $\Omega\cap\{t>\psi(x)\}$,   
 \item the restriction of $v$ on 
 $\Sigma$ coincides with $v_0$.
 \end{itemize}
 Such a function $v(t, x)$ is unique in a neighborhood of $\Sigma$. 
\end{theorem}

\begin{remark}
\noindent
[A] \  
If $f_2(\tau, \xi)= a^{-1} (\tau^2-|\xi|^2)$, then 
  (\ref{eq:assumption}) reduces to 
$1-|\nabla\psi(x)|^2\ne 0$. 
It is satisfied by almost any function $\psi$. 

\noindent
[B] \ 
The function $\psi$ is a solution 
not to the Eikonal equation but to  
another  equation 
of the first order. 
It would be an interesting task 
to formulate a  general theory encompassing 
the "pseudo-Eikonal" equations 
(\ref{eq:assumption}), (\ref{eq:assumption_higher_1}) 
and the curvature condition in \cite{KL1}. 
A promising approach is   to  deal  with them in the  
framework of \cite{KS}, but it is beyond the scope 
of the present paper.

\noindent
[C] \
The speed of blowup is controlled by the 
coefficient $a$. It is calculated in terms of 
$f_2$ and $\psi$. 
On the other hand,   we can find many 
$\psi$'s satisfying (\ref{eq:assumption}) 
for a given pair $(f_2, a)$. 
We have only to use  the Cauchy-Kovalevskaya theorem or the 
classical theory of  partial 
differential equations of the first order. 

Recall that the speed of blowup was always the 
same irrespective of   $\psi$ 
in the study of $\square u=e^u$ in \cite{KL1} and \cite{KL2}.  
\end{remark}

\section{Proof of Theorem~\ref{thm:main}}
\label{sec:proof}

We introduce a new system of coordinates $(T, X)$ by 
$T=t-\psi(x),\, X_i=x_i\,(i=1, 2, \dots, n)$. 
The derivatives in  the new variables are  marked with   
$ \hat {}$ :
\[
 \partial_i=\partial/\partial x_i, \;
 \partial_t=\partial/\partial t, \;
 {\hat {\partial}}_i=\partial/\partial X_i, \;
 {\hat {\partial}}_X=({\hat {\partial}}_1, \dots, {\hat {\partial}}_n), \;
 {\hat {\partial}}_T=\partial/\partial T. 
\]
In the same way, $\Delta$ and $\nabla$ are the 
Laplacian and the gradient in $x$ respectively, 
while $ \widehat \Delta$ and $ \widehat \nabla$ are 
those in $X$. 
Subscripts on functions refer to  the derivatives in $x$. 
For example, $\psi_i=\partial_i \psi, \psi_{ii}=\partial_i^2 \psi$. 

The derivatives transform according to the following rules:
\[
 \partial_t={\hat {\partial}}_T, 
 \quad
 \partial_i=-\psi_i {\hat {\partial}}_T + {\hat {\partial}}_i . 
\]

We want to express the d'Alembertian $\square=\partial_t^2-\Delta$ 
in terms of ${\hat {\partial}}_i$ and  ${\hat {\partial}}_T$. 
This calculation has already been explained in \cite{KL1}. 
We nevertheless repeat the argument here for the sake of 
the reader. It might seem a little confusing since 
it employs two systems of coordinates simultaneously. 
To begin with, we get
\begin{align*}
 \partial_i^2
 & = 
 \partial_i (-\psi_i {\hat {\partial}}_T+{\hat {\partial}}_i) 
 =-\partial_i (\psi_i {\hat {\partial}}_T) + \partial_i{\hat {\partial}}_i
 \\
 & =-\psi_i \partial_i {\hat {\partial}}_T - \psi_{ii}{\hat {\partial}}_T 
 + (-\psi_i {\hat {\partial}}_T + {\hat {\partial}}_i) {\hat {\partial}}_i 
  \\
 &=-\psi_i (-\psi_i {\hat {\partial}}_T + {\hat {\partial}}_i) {\hat {\partial}}_T
 -\psi_{ii}{\hat {\partial}}_T -\psi_i{\hat {\partial}}_T {\hat {\partial}}_i 
 +{\hat {\partial}}_i^2
  \\
 & =
 \psi_i^2 {\hat {\partial}}_T^2
 - 2\psi_i {\hat {\partial}}_i {\hat {\partial}}_T 
 -\psi_{ii} {\hat {\partial}}_T 
 +{\hat {\partial}}_i^2 
  .
\end{align*}
Therefore $\square$ has the following form:
\begin{equation}
\label{eq:d'Alembert}
 \square
 = \Psi {\hat {\partial}}_T^2 
 +2\sum_{i=1}^n \psi_i {\hat {\partial}}_i{\hat {\partial}}_T  + (\Delta \psi){\hat {\partial}}_T
 - \widehat \Delta, 
\end{equation}
where $\Psi =\Psi(X)=\Psi(x)=  1-|\nabla \psi(x)|^2\ne 0$.

Set $u(t, x)=- a \log  T+v$, where 
$v$ is an unknown analytic function near 
$\Sigma=\{T=0\}$. 
Then by (\ref{eq:d'Alembert}), the left hand side of 
(\ref{eq:wave}) is 
\begin{equation}
  \square u
 =
 \frac{a \Psi}{T^2} - \frac{a \Delta\psi}{T}+\square v. 
 \label{eq:lhs}
\end{equation}
Next, let us calculate the right hand side. 
We have 
\[
 \partial_t u =-a/T + {\hat {\partial}}_T v, \quad
 \partial_i u = a \psi_i /T -\psi_i {\hat {\partial}}_T v + {\hat {\partial}}_i v.
\] 
It leads to  
\begin{align}
  (\partial_t u)^2
 &=\frac{a^2}{T^2}
  -\frac{2a}{T}{\hat {\partial}}_T v + ({\hat {\partial}}_T v)^2, 
 \nonumber
 \\
 (\partial_t u)(\partial_i u)
 &= - \frac{a^2 \psi_i}{T^2} 
 +\frac{2a\psi_i}{T}{\hat {\partial}}_T v
 -\frac{a}{T}{\hat {\partial}}_i v 
 +({\hat {\partial}}_T v) (-\psi_i {\hat {\partial}}_T v+{\hat {\partial}}_i v), 
  \nonumber
  \\
 (\partial_i u) (\partial_j u) 
 &=
 \frac{a^2 \psi_i \psi_j}{T^2}
 +\frac{1}{T}
 \left\{
   -2a\psi_i \psi_j {\hat {\partial}}_T v 
   +a (\psi_i {\hat {\partial}}_j + \psi_j {\hat {\partial}}_i)v
 \right\}
 \nonumber
 \\
 &\qquad
 +(-\psi_i {\hat {\partial}}_T v + {\hat {\partial}}_i v)
  (-\psi_j {\hat {\partial}}_T v + {\hat {\partial}}_j v).
  \nonumber
\end{align}
On the other hand, 
there exists an analytic function $\tilde{f_2}(T, X)$  
in a neighborhood of $\Sigma$ such that 
\[
  f_2 (t, x; -1, \nabla \psi(x))
  =  f_2 (\psi(x), x; -1, \nabla \psi(x))
  +T \tilde{f_2}(T, X). 
\]
We set 
$f_2 (\Sigma)=f_2 (\psi(x), x; -1, \nabla \psi(x))$ for brevity. 
This quantity appears in (\ref{eq:assumption}) and 
should be distinguished from 
$f_2 (-1, \nabla \psi(x))= f_2 (t, x; -1, \nabla \psi(x))$, 
a particular case of 
$f_l (\tau, \xi)= f_l (t, x; \tau, \xi)$. 
We get 
\[
  \frac{f_2 (-1, \nabla \psi(x))}{T^2}
  =\frac{f_2(\Sigma)}{T^2}+\frac{\tilde{f_2}(T, X)}{T}. 
\]
Consequently,  the right hand side of  
the equation (\ref{eq:wave}) is 
\begin{align}
 f (t, x; \partial_t u, \nabla u)
 &=
 \frac{a^2 f_2 (\Sigma)}{T^2}
 + \frac{a^2 \tilde{f_2}(T, X)}{T}
 -\frac{2a}{T}f_2 (-1, \nabla\psi){\hat {\partial}}_T v
 \label{eq:rhs}
 \\
 &\quad 
 +\frac{a}{T}L(T, X, {\hat {\partial}}_X)v
 +f_2 
 \left( {\hat {\partial}}_T v, \{-\psi_i {\hat {\partial}}_T v+{\hat {\partial}}_i v\}_i
 \right)
 \nonumber
 \\
 &\quad 
 + \frac{a}{T} f_1 (-1, \nabla\psi)
 +f_1 \left(
           {\hat {\partial}}_T v, \{-\psi_i {\hat {\partial}}_T v+{\hat {\partial}}_i v\}_i
      \right)
 +f_0 (t, x)
 , 
 \nonumber
\end{align}
where $L(T, X, {\hat {\partial}}_X)$ is a vector field 
with $[L, T]=0$. 

Now let us compare (\ref{eq:lhs}) and (\ref{eq:rhs}). 
Our assumption (\ref{eq:assumption}) implies that 
 the coefficients of $T^{-2}$ coincide: $a\Psi=a^2 f_2 (\Sigma)$. 
Therefore, in view of (\ref{eq:d'Alembert}),  
the equation (\ref{eq:wave}) is equivalent to 
the following one:
\begin{align}
 &
 - a \Delta\psi
 + T 
     \left\{
         \Psi{\hat {\partial}}_T^2 v 
         +2\sum_{i=1}^n \psi_i {\hat {\partial}}_i {\hat {\partial}}_T  v
         + (\Delta \psi){\hat {\partial}}_T v  - \widehat \Delta v 
     \right\}
 \label{eq:Fuchs} 
 \\
 &=
 a^2 \tilde{f}_2 (T, X)
 - 2\Psi {\hat {\partial}}_T v
  + a  L(T, X, {\hat {\partial}}_X)v
 +T f_2 
 \left( {\hat {\partial}}_T v, \{-\psi_i {\hat {\partial}}_T v+{\hat {\partial}}_i v\}_i
 \right)
 \nonumber
 \\
 &\quad 
 +  a f_1 (-1, \nabla\psi)
 +T f_1 \left(
           {\hat {\partial}}_T v, \{-\psi_i {\hat {\partial}}_T v+{\hat {\partial}}_i v\}_i
      \right)
 +T f_0 (t, x) .
 \nonumber
\end{align}
We divide this equation by $\Psi\ne 0$ 
and get an example of 
(\ref{eq:GT_weight1}) in \S\ref{sec:nonlinFuchs} below.  
By Proposition~\ref{pro:GT}, 
there exists a unique 
analytic solution  $v$ to (\ref{eq:Fuchs}) with 
$v\bigr|\raisebox{-0.7ex}{}_\Sigma  = v_0$ in a neighborhood of  $\Sigma$. 
The function  
$u=- a \log  T + v$ solves (\ref{eq:wave}) in $\{T>0\}$. 
We have finally proved Theorem~\ref{thm:main}.

\section{Nonlinear Fuchsian partial differential equations}
\label{sec:nonlinFuchs}

In this section, we   prove a result on 
nonlinear  Fuchsian partial differential 
equations which was used in the proof of  Theorem~\ref{thm:main}.

We consider a nonlinear equation 
in a neighborhood of 
\( 
       \{0\}\times V \subset 
       \mathbb{R}_t \times \mathbb{R}^n_x
\), 
where $V$ is an open set of $\mathbb{R}^n_x$. 
We employ the usual notation for derivatives:  
\( \partial_i=\partial/\partial x_i, \partial_x = (\partial_1, \dots, \partial_n), 
  \nabla v=(\partial_1 v, \dots, \partial_n v)
\) 
and 
$\partial_t=\partial/\partial t$. 

Now let us consider a Cauchy problem for the 
following  nonlinear Fuchsian partial differential equation of 
weight 1:
\begin{align}
 t \partial_t^2 v + 2 \partial_t v
 &=
 \alpha (t, x)
 + t P (t, x,  \partial_x) \partial_t v + Q(t, x,  \partial_x) v 
 \label{eq:GT_weight1}
 \\
 &\qquad
 + t \beta(t, x) (\partial_t v)^2
 + R(t, x, \partial_x) v \cdot t \partial_t v 
 +  S(t, x; \nabla v).
 \nonumber
\end{align}
Here  $\alpha$ and $\beta$ are  analytic functions. 
The linear differential operators 
$P, Q$ and  $R$  with analytic coefficients are  
of order $1, 2$ and  $1$ respectively 
and 
$[P, t]=[Q, t]=[R, t]=0$. 
We assume that 
$S(t, x; \xi)$ is a quadratic form in $\xi$ 
whose  coefficients are analytic functions in $(t, x)$.

\begin{proposition} 
Let $v_0 (x)$ be an   analytic function in  $V$. 
Then 
(\ref{eq:GT_weight1}) has a unique analytic solution 
$v(t, x)$ with $v(0, x)=v_0 (x)$ 
near $\{0\}\times V \subset\mathbb{R}\times\mathbb{R}^n$.
\label{pro:GT}
\end{proposition}

\begin{proof}
First we reduce the problem to the case $v_0 (x)\equiv 0$. 
We introduce a new unknown function $w$ by 
 $v(t, x)=v_0 (x)+w(t, x)$. Then we have 
$t \partial_t^2 v + 2 \partial_t v=t \partial_t^2 w + 2 \partial_t w$. 

We denote the right hand side  of 
(\ref{eq:GT_weight1}) by $N[v]$. 
An elementary computation shows 
\begin{align*}
 N[v]&=N[v_0 + w]
 \\
 &
 =N[w]
  +Q v_0 + R v_0 \cdot t \partial_t w 
  + S (t, x; \nabla v_0)  + T(t, x; \nabla w), 
\end{align*}
where $T(t, x; \xi)=\sum_{i=1}^n c_i(t, x)\xi_i$ for some 
analytic function $c_i (t, x)$.   
This implies that $w(t, x)$ satisfies an equation 
of the type (\ref{eq:GT_weight1}), possibly 
with different coefficients. 
Therefore we have only to prove that 
(\ref{eq:GT_weight1}) has a unique analytic solution 
$v(t, x)$ with $v(0, x)\equiv 0$. 
Multiplication by $t$ yields a   nonlinear Fuchsian 
partial differential equation solved in 
\cite{GT_art} and 
\cite{GT_book}. 
Its characteristic polynomial is $\rho^2+\rho$ and 
there exists a unique analytic solution 
$v(t, x)$ with $v(0, x)\equiv 0$. 
%
\end{proof}

\begin{remark}

\noindent
[A] \  
There is a stronger uniqueness result due to 
G\'erard-Tahara. We only give a sketch here 
and refer the reader to   \cite{GT_art} or 
\cite{GT_book} for 
a precise statement. 

Assume that $v_0 (x) \equiv 0$ for simplicity 
and let $\tilde v(t, x)$ be a multi-valued holomorphic solution 
to (\ref{eq:GT_weight1}) in the universal covering space of 
$\{t\ne0\}\subset \mathbb{C}_t \times \mathbb{C}_x^n$. 
If there exists a positive constant $p$ such that 
$\tilde v(t, x)=O(|t|^p)$ as $t$ tends to zero, then 
$\tilde v(t, x)$ is nothing but the solution $v(t, x)$ in 
Proposition~\ref{pro:GT}. 

This observation enables us to improve the uniqueness part 
of Theorem~\ref{thm:main}. It rules out some functions 
involving terms of the form $T^\alpha \log^\beta T, \alpha>0, T=t-\psi(x)$.

\noindent
[B] \ The present  formulation, a very restrictive one,  
has been chosen  because 
it is all we  need in the proof of 
Theorem~\ref{thm:main}.  
There is much room for generalization.   
We can easily consider other  characteristic polynomials than  
$\rho^2+\rho$. 
Difficulties arise 
if we consider  nonlinear terms 
which are not entire in its arguments; 
it is not always possible to consider 
large initial values. 
\end{remark}

\section{Time-reversal}
\label{sec:time-reversal}

In this section we discuss time-reversal and construct 
solutions on the negative side of $\Sigma$ under 
certain conditions.  
The equation (\ref{eq:wave}) is not invariant under 
time-reversal. Indeed, if we set $s=-t$, it becomes 
\begin{equation}
\label{eq:timereversal1}
 u_{ss}-\Delta u= f(-s, x, -\partial_s u, \nabla u).
\end{equation}
In addition to (\ref{eq:assumption}), we assume 
\begin{equation}
\label{eq:timereversal2}
 f_2 (t, x; 1, \nabla\psi)=f_2 (t, x; -1, \nabla\psi).
\end{equation}
Here notice that $f_2 (\tau, \xi)=f_2 (-\tau, -\xi)$ is trivial. 
The condition (\ref{eq:timereversal2}) is satisfied 
if $f_2 (\tau, \xi)$ is free of terms of the form 
$c_j (t, x)\tau\xi_j$.

We set 
\(
 \tilde\psi(x)=-\psi(x)
\) and introduce 
\( 
   \tilde f_2 (\tau, \xi)= f_2 (-\tau, \xi)
\), which is, to (\ref{eq:timereversal1}),  
what $f_2$ is to (\ref{eq:wave}). 
Combination of (\ref{eq:assumption}) and 
(\ref{eq:timereversal2}) yields 
\begin{equation}
\label{eq:assumptiontimereversal}
 0\neq
  1-|\nabla\tilde\psi(x)|^2
  =
  a \tilde f_2 \bigl(-1, \nabla\tilde\psi(x)\bigr).
\end{equation}
By using Theorem~\ref{thm:main}, we can find   solutions 
to (\ref{eq:timereversal1}) in 
$\{(s, x); s>\tilde\psi(x)\}$. 
It gives a solution to (\ref{eq:wave}) in 
$\{(t, x); t< \psi(x)\}$, the negative side of $\Sigma$. 

To sum up, 
under (\ref{eq:assumption}) and 
(\ref{eq:timereversal2}), there exist    solutions defined 
on both sides of   
the  noncharacteristic  hypersurface  $\Sigma$ 
and their  singularities propagate along it.

\section{Elliptic case}
\label{sec:elliptic}

We can apply a similar argument to 
\begin{equation}
\label{eq:elliptic}
  \Delta u=a^{-1} |\nabla u|^2
\end{equation}
in $\mathbb{R}^n$, where $a$ is a nonzero constant. 
In this case, 
we can treat any analytic nonsingular hypersurface $\Lambda$ 
at least locally. 
Indeed, after rotation, it can be written in the 
form $\Lambda=\{x_1=\phi(x_2, \dots, x_n)\}$. 
Notice that the equation is rotation-invariant. 
The function $1+\sum_{i=2}^n (\partial \phi/\partial x_i)^2>0$ 
replaces $\Psi$. 
Hence we have the following:
\begin{theorem}
 Let $p$ be a point of an arbitrary 
 analytic nonsingular hypersurface $\Lambda$. 
 Then there exists a neighborhood  $U$ of $p$ such that 
 in a connected component of $U\setminus \Lambda$
 there exists an analytic solution $u(x)$ to (\ref{eq:elliptic}) 
 which blows up along $\Lambda\cap U$. 
\end{theorem}

\section{Multi-valued solutions in the complex domain}
\label{sec:Multi-valued}

Our calculus is valid in the complex domain 
and the solution $u(t, x)$ in Theorem~\ref{thm:main} 
becomes a multi-valued holomorphic  function around the 
\textit{noncharacteristic} 
complex hypersurface 
$\Sigma^\mathbb{C}=\{(t, x)\in\mathbb{C}^{n+1};     t=\psi(x)\}$. 
It is the unique solution to the following  singular Cauchy problem 
near a point $(t_0, p)\in\Sigma^\mathbb{C}, t_0=\psi(p)$: 
\[
 \begin{cases}
  &\square u=f(t, x, \partial_t u, \nabla u), 
  \\
  &
  u(t_0, x)=-a \log (t_0 - \psi (x))+ v (t_0, x), 
  \\
  & \partial_t u(t_0, x)=-a/(t_0 - \psi (x))+ \partial_t v(t_0, x). 
 \end{cases}
\]
For such Cauchy data, 
singularities propagate along the noncharacteristic hypersurface 
$\Sigma^\mathbb{C}$ and only along it.  
Notice that some results in \cite{KL1} and \cite{KL2} 
are open to  similar interpretations. 

Several authors have constructed multi-valued solutions to 
nonlinear equations around \textit{characteristic} hypersurfaces. 
See  \cite{GT_art},  \cite{GT_book}, \cite{KL2} and 
other works 
by Kichenassamy and by Tahara about  Fuchsian problems. 
They give   solutions which are singular along  the initial surface 
$\{t=0\}$. 
Moreover, see \cite{Uchikoshi} and the references cited therein for 
singularities obtained by solving 
Eikonal equations.

\section{Equations with nonlinearity of higher order}

The leading term of a singular solution can be guessed by 
solving an ordinary differential equation. 
For example, one finds a solution 
$u(t)=\log 2-2\log t$ to $d^2 u / dt^2 = e^t$ 
and it was the basis of \cite{KL1} and \cite{KL2}. 
On the other hand, the solution 
$u(t)=-\log t$ to 
$d^2 u / dt^2 = (du/dt)^2$ is the simplest case of 
our Theorem~\ref{thm:main}.

We are naturally tempted to solve 
$d^2 u / dt^2 = -(du/dt)^{m+1}$ for $m\ge 2$ 
and get 
$u(t)=C_m t^{(m-1)/m}$, where $C_m=m^{(m-1)/m}/(m-1)$. 
This solution vanishes as $t>0$ tends to $0$, but 
its derivative   blows up. 
In this section, we give a multi-dimensional version of this fact.

Let $f(t, x; \tau, \xi)$ be a polynomial of 
degree $m+1\, (m\ge 2)$ in 
$(\tau, \xi)\in\mathbb{R}\times\mathbb{R}^n$  whose 
coefficients are analytic functions in 
$(t, x)$ in an open set $U$ of $\mathbb{R}\times\mathbb{R}^n$. 
We decompose $f$ into homogeneous parts: 
$f=\sum_{l=0}^{m+1}f_l$, where $f_l$ is  homogeneous of 
degree $l$. 
We sometimes denote $f_l (t, x; \tau, \xi)$ by 
$f_l (\tau, \xi)$ for short. 

We consider 
\begin{equation}
 u_{tt}-\Delta u=f(t, x; \partial_t u, \nabla u)
 =
 \sum_{l=0}^{m+1}f_l (t, x; \partial_t u, \nabla u) 
\label{eq:wave_higher}
\end{equation}
and  construct a solution 
$u(t, x)$   which  vanishes on a prescribed 
\textit{noncharacteristic} hypersurface. 

\begin{theorem}
\label{thm:higher}
 Let $\Sigma$ be a hypersurface defined by 
 \(\Sigma=\{ 
             (t, x)\in U;     t=\psi(x)
          \}
 \). 
 Assume   
 \begin{align}
    &0\ne 1-|\nabla\psi(x)|^2 
    = 
      \frac{(-m+1)^m a^m}{m^{m-1}}
     f_{m+1} \bigl(\psi(x), x; -1,  \nabla \psi(x)\bigr), 
    \label{eq:assumption_higher_1}
    \\
    & f_m \bigl(\psi(x), x ;
       -1,  \nabla\psi (x)\bigr)
           =0
    \label{eq:assumption_higher_2}
 \end{align}
 for  a nonzero constant $a$.  
 Then 
 there exists an open neighborhood $\Omega$ of 
 $\Sigma$ and an  analytic solution  $u(t, x)$ 
 to (\ref{eq:wave_higher})  
 in  $\Omega\cap\{t>\psi (x)\}$ 
    such that: 
 \begin{align*}
     u(t, x)
     &
     \sim a \bigl(t-\psi(x)\bigr)^{(m-1)/m}
     \ \mbox{(vanishing)}, 
     \\
     \partial_t u(t, x)
     &
     \sim a 
     \frac{m-1}{m}
     \bigl(t-\psi(x)\bigr)^{-1/m}
     \ \mbox{(blow-up)}
 \end{align*}
  as $t-\psi(x)\to +0$. 
  This solution $u(t, x)$ is analytic in 
  $\bigl((t-\psi)^{1/m}, x\bigr)$
  and no logarithm is involved. 
\end{theorem}

\begin{remark}
The conditions (\ref{eq:assumption_higher_1}) 
and (\ref{eq:assumption_higher_2}) are 
equivalent to the existence of such 
analytic functions $\tilde{f}_{m+1}(T, X)$ and $\tilde{f}_{m}(T, X)$  
that 
\begin{align*}
  &f_{m+1} \bigl(t, x ;
       1, -\nabla\psi (x)\bigr)
  =
  \frac{m^{m-1}}{(-m+1)^m a^m} 
  (1-|\nabla\psi(x)|^2)
  +T \tilde{f}_{m+1}(T, X)
  , 
  \tag{\ref{eq:assumption_higher_1}'}   
  \\
  &f_m \bigl(t, x ;
       1, -\nabla\psi (x)\bigr)
  = T  \tilde{f}_m (T, X)
  \tag{\ref{eq:assumption_higher_2}'}.    
\end{align*}
\end{remark}

\smallskip

\begin{proof}
We introduce a system of coordinates $(T, X)$ as before 
and seek a solution $u$ of the form 
\begin{equation}
 u=a T^{(m-1)/m} + T v, 
\end{equation}
where $v$ is an unknown  function 
analytic in $(T^{1/m}, X)$. 
It may have fractional-power singularities as a 
function in $T$. 
By using the formulas 
\[
 {\hat {\partial}}_T T = T {\hat {\partial}}_T +1, \quad
 {\hat {\partial}}_T^2 T = T {\hat {\partial}}_T^2 + 2{\hat {\partial}}_T
\]
and (\ref{eq:d'Alembert}), 
we can express the 
left hand side of (\ref{eq:wave_higher}) as follows:
\begin{align}
 \square u
  =&  
   L_1
    + a \frac{m-1}{m} (\Delta \psi) T^{-1/m}
  \label{eq:lhs_higher}
 \\
 & +  
 \Bigl\{
      \Psi (T {\hat {\partial}}_T^2 + 2{\hat {\partial}}_T)
      +2 \sum_{i=1}^n \psi_i (T {\hat {\partial}}_T {\hat {\partial}}_i + {\hat {\partial}}_i) 
      +(\Delta\psi) (T {\hat {\partial}}_T + 1) 
      - T \widehat\Delta
    \Bigr\} v, 
  \nonumber
\end{align}
where   
\(
 L_1=a \dfrac{-(m-1)}{m^2} \Psi T^{-(m+1)/m}
\).

Next we shall compute the right hand side of (\ref{eq:wave_higher}). 
We get 
\[
\begin{cases}
 &\partial_t u 
  = 
  a \frac{m-1}{m} T^{-1/m} + (T {\hat {\partial}}_T v +v), 
 \\
 &\partial_i u 
  = -\psi_i \partial_t u+ T {\hat {\partial}}_i v.
\end{cases}
\]
Notice that the expression of $\partial_i u$ contains 
$\partial_t u $. 
The latter quantity   plays an important role in the 
following calculation. 
We have
\begin{align}
\label{eq:f_l}
 f_l (\partial_t u, \nabla u)
 &= f_l 
   \bigl(\partial_t u, 
         \{
            -\psi_i \partial_t u + T {\hat {\partial}}_i v
         \}_i
   \bigr)
 \\
 \nonumber
 &= 
 f_l (1, -\nabla \psi) (\partial_t u)^l 
 + R_l, 
\end{align}
where $R_l$ is a homogeneous polynomial of 
degree $l$ in $\partial_t u$ and $T {\hat {\partial}}_i v$ and 
is of degree $\le \max\{0, l-1\} (\le m)$ with respect to $\partial_t u$. 
It is 
of  degree $\ge 1$ with respect to $T {\hat {\partial}}_i v$ 
if $l\ge 1$. 
The "spare" factor $T$ in $T {\hat {\partial}}_i v$ helps us since 
we are trying to get a Fuchsian equation; 
it offsets the negative power of 
$T$ in $\partial_t u$. 

We have
\begin{equation}
\label{eq:dtul}
 (\partial_t u)^l 
   = 
 \left(a \frac{m-1}{m} \right)^l 
 T^{-l/m}
 + 
 l 
 \left(a \frac{m-1}{m} \right)^{l-1} 
 T^{-(l-1)/m}
 (T {\hat {\partial}}_T v + v) 
 +B_l, 
\end{equation} 
where $B_l$ is the remainder of the binomial expansion, that is, 
\[
 B_l = 
\sum_{j=0}^{l-2} 
 \binom{l}{j}
 \left( a\frac{m-1}{m} \right)^{j}
 T^{-j/m} 
 (T {\hat {\partial}}_T v + v)^{l-j}. 
\]
Since 
$f=\sum_{l=0}^{m+1}f_l$, we have 
by (\ref{eq:f_l}) and (\ref{eq:dtul}) 
\begin{align}
 &f(\partial_t u, \nabla u)
 \label{eq:rhs_higher}
  \\
 = & 
 \left\{
  L_2 + L_3  
  + f_{m+1}  (1, -\nabla \psi) B_{m+1}
 \right\}
  +\sum_{l=0}^{m+1} R_l
  +\sum_{l=0}^{m} f_l (1, -\nabla \psi) (\partial_t u)^l , 
  \nonumber
\end{align}
where 
\begin{align*}
 L_2
 &=
  f_{m+1} (1, -\nabla \psi(x))
  \left(a \frac{m-1}{m} \right)^{m+1}
  T^{-(m+1)/m},   
 \\
 L_3 
 &=
  f_{m+1}  (1, -\nabla \psi(x)) \cdot
  (m+1)  \left( a \frac{m-1}{m} \right)^m 
  T^{-1} (T{\hat {\partial}}_T v + v) .
\end{align*}
We set 
\[
  L'_2 =f_{m+1}(\psi(x), x: 1, -\nabla\psi(x))
        \left(a \frac{m-1}{m} \right)^{m+1} 
        T^{-(m+1)/m}  
\]
then there exists  an analytic function $L''_2 (T, X)$ 
such that 
\[
   L_2=L'_2 + T^{-1/m} L''_2 (T, X). 
\]
The assumption (\ref{eq:assumption_higher_1}) 
implies that 
$L_1=L'_2$ and that 
\[
 L_3
 = 
 -\Psi\cdot\frac{m+1}{m}T^{-1} 
 (T {\hat {\partial}}_T v + v).
\]
Note that (\ref{eq:assumption_higher_1}) is 
written in terms of $(-1, \nabla\psi)$ instead of 
$(1, -\nabla\psi)$. 
By subtracting $L_1=L'_2$ from 
(\ref{eq:wave_higher}), we come 
very close to a Fuchsian equation 
in view of (\ref{eq:lhs_higher}) and (\ref{eq:rhs_higher}).

We set $s =T^{1/m}$, ${\hat {\partial}}_s =\partial/\partial s$ in $T>0$. 
Since $m T {\hat {\partial}}_T = s {\hat {\partial}}_s $, 
we have 
\begin{align*}
  T{\hat {\partial}}_T v + v
  &=
  \frac{1}{m}   (s {\hat {\partial}}_s v+mv), 
  \\
  \partial_t u
  &=
  \frac{1}{m} s^{-1}
   \{
    a(m-1)+s ( s {\hat {\partial}}_s v + mv)
   \}.
\end{align*}
Now we   multiply the equation 
$\square u -L_1=f-L'_2$ by $m^2 s^m/\Psi$ 
and   obtain a  nonlinear Fuchsian equation which 
belongs to the class treated  
in \cite{GT_art} and \cite{GT_book}. 
Namely it has the form:
\begin{align*}
 &  \frac{a m(m-1)}{\Psi}(\Delta\psi)s ^{m-1}
 +
 \bigl\{(s {\hat {\partial}}_s )^2 + ms {\hat {\partial}}_s \bigr\}
 v
 \\
 & \qquad 
  + s ^m P_1 
 \left(s,  X; 
               \{(s {\hat {\partial}}_s +m){\hat {\partial}}_i  v\}_i, 
                 s {\hat {\partial}}_s v +mv, \widehat\Delta v
 \right) 
 \\
 = 
 &
  s^{m-1} \tilde{L}_2 (s, X)
  -(m+1)(s {\hat {\partial}}_s v +mv)
  +
  s P_2 (s,  X;  s {\hat {\partial}}_s v +mv  )
  (s {\hat {\partial}}_s v + mv)^2        
  \\
  &\qquad
  + s ^m P_3 \left(s,  X;  s {\hat {\partial}}_s v   + mv, 
                \{{\hat {\partial}}_i v\}_i  
               \right)     
  \\
  &\qquad
  +\sum_{l=0}^m Q_l (s,  X) 
  s ^{\max\{m-l, 1\}}
  \bigl\{
     a(m-1)+s ( s {\hat {\partial}}_s v + mv)
  \bigr\}^l
  .
\end{align*}
Here $\tilde{L}_2 (s, X)$ is 
an analytic function and 
$P_1, P_2$ and $P_3$ are polynomials in the 
arguments coming after the semicolons. 
Their coefficients and $Q_l$ are analytic functions in $(s,  X)$ near 
$s =0$. 
The exponent of $s $ immediately after $Q_l$ is not   $m-l$ but 
$\max\{m-l, 1\}$ because of (\ref{eq:assumption_higher_2}). 
The characteristic polynomial is $(\rho+m)(\rho+m+1)$ 
and 
there exists a solution $v$, which is an 
analytic function in $(s, X)=((t-\psi)^{1/m}, x)$. 
\end{proof}

\end{document}